\documentclass[reqno,12pt]{article}%
\usepackage{amsmath}
\usepackage{graphicx}
\usepackage{amsfonts}
\usepackage{amssymb}
\usepackage{latexsym}%
 \makeatletter
\newtheorem{theorem}{Theorem}

\newtheorem{claim}[theorem]{Claim}

\newtheorem{example}[theorem]{Example}

\newtheorem{lemma}[theorem]{Lemma}

\newenvironment{proof}[1][Proof]{\noindent{\textbf {#1}  }}  {\hfill$\Box$\bigskip}

\begin{document}

\title{Chromatic number and minimum degree of $K_{r}$-free graphs}
\author{Vladimir Nikiforov\\{\small Department of Mathematical Sciences, University of Memphis, Memphis,
TN 38152}\\{\small email: vnikifrv@memphis.edu}}
\maketitle

\begin{abstract}
Let $\delta\left(  G\right)  $ be the minimum degree of a graph $G.$ A number
of famous results about triangle-free graphs determine the maximum chromatic
number of graphs of order $n$ with $\delta\left(  G\right)  >n/3.$ In this
note these results are extended to $K_{r+1}$-free graphs of order $n$ with
$\delta\left(  G\right)  >\left(  1-2/(2r-1\right)  n.$ In particular:

(a) there exist $K_{r+1}$-free graphs of order $n$ with $\delta\left(
G\right)  >\left(  1-2/\left(  2r-1\right)  \right)  n-o\left(  n\right)  $
and arbitrary large chromatic number;

(b) if $G$ is a $K_{r+1}$-free graph of order $n$ with $\delta\left(
G\right)  >\left(  1-2/\left(  2r-1\right)  \right)  n,$ then $\chi\left(
G\right)  \leq r+2;$

(c) the structure of the $\left(  r+1\right)  $-chromatic $K_{r+1}$-free
graphs of order $n,$ with $\delta\left(  G\right)  >\left(  1-2/(2r-1\right)
n$ is found.\medskip


\end{abstract}

\section{Introduction and main results}

In notation we follow \cite{Bol98}. In 1962 Andr\'{a}sfai \cite{And62}
introduced the function
\[
\psi\left(  n,r,h\right)  =\max\left\{  \delta\left(  G\right)  :G\text{ is a
}K_{r+1}\text{-free graph of order }n\text{ with }\chi\left(  G\right)  \geq
h\right\}  ,
\]
which has been widely studied during the years. One of the first contributions
were the famous theorem and example of Andr\'{a}sfai, Erd\H{o}s and S\'{o}s
\cite{AES74} showing that for every $r\geq2,$
\begin{equation}
\left(  1-\frac{3}{3r-1}\right)  n+O\left(  1\right)  \leq\psi\left(
n,r,r\right)  \leq\left(  1-\frac{3}{3r-1}\right)  n. \label{AESbas}%
\end{equation}

Another milestone along this road is an example of Hajnal \cite{ErSi73}
showing that $\psi\left(  n,2,h\right)  >n/3-o\left(  n\right)  $ for every
$h\geq3;$ for an updated version of Hajnal's example see Example \ref{exEHS}
below. Thus, for $r=2,$ this example leaves unanswered only one simple, yet
tricky question: \emph{how large can be }$\chi\left(  G\right)  $\emph{ of a
}$K_{3}$\emph{-}$f\emph{ree}$ \emph{graph }$G$\emph{ of order }$n$
$\emph{with}$ \emph{ }$\delta\left(  G\right)  >n/3.$\emph{ }

Erd\H{o}s and Simonovits \cite{ErSi73} conjectured that all $K_{3}$-free
graphs of order $n$ with $\delta\left(  G\right)  >n/3$ are $3$-chromatic, but
this was disproved by H\"{a}ggkvist \cite{Hag82}, who described for every
$k\geq1$ a $10k$-regular, $4$-chromatic, $K_{3}$-free graph of order $29k,$
for a description see Example \ref{exHag} below. The example of H\"{a}ggkvist
is based on the Mycielski graph $M_{3},$ also known as the Gr\"{o}tzsch graph,
which is a $4$-chromatic $K_{3}$-free graph of order $11$. We shall see later
that $M_{3}$ is a true landmark in this area, but let us first recall the
Mycielski graphs given in \cite{Myc55}: a sequence $M_{1},M_{2},\ldots$ of
$K_{3}$-free graphs with $\chi\left(  M_{i}\right)  =i+1,$ constructed as
follows: \medskip

\emph{Set }$M_{1}=K_{2}.$\emph{ To obtain }$M_{i+1}:$\emph{ write }%
$v_{1},\ldots,v_{n}$\emph{ for the vertices of }$M_{i};$\emph{ choose }%
$n+1$\emph{ other vertices }$u_{1},\ldots,u_{n+1};$\emph{ for every }%
$i\in\left[  n\right]  $\emph{ join }$u_{i}$\emph{ precisely to the neighbors
of }$v_{i};$\emph{ join }$u_{2n+1}$\emph{ to }$u_{1},\ldots,u_{n}.$\medskip

Other graphs crucial in the study of $\psi\left(  n,r,h\right)  $ are the
$K_{3}$-free $3$-chromatic Andr\'{a}sfai graphs $A_{1},A_{2},\ldots,$ first
described in \cite{And62}:\medskip\ 

\emph{Set }$A_{1}=K_{2}$\emph{ and for every }$i\geq2$\emph{ let }$A_{i}%
$\emph{ be the complement of the }$\left(  i-1\right)  $\emph{'th power of
}$C_{3i-1}.$ \medskip

In \cite{Jin93}, Jin proved the following theorem, generalizing the case $r=2$
of the theorem of Andr\'{a}sfai, Erd\H{o}s and S\'{o}s and a result of
H\"{a}ggkvist from \cite{Hag82}.\medskip

\textbf{Theorem A }\emph{Let }$1\leq k\leq9,$\emph{ and let }$G$\emph{ be a
}$K_{3}$\emph{-free graph of order }$n.\emph{\ }$\emph{If}
\[
\delta\left(  G\right)  >\frac{k+1}{3k+2}n,
\]
\emph{ then }$G$\emph{ is homomorphic to }$A_{k}.$\medskip

Note that this result is tight: taking the graph $A_{k+1}$, and blowing it up
by a factor $t,$ we obtain a $K_{3}$-free graph $G$ of order $n=\left(
3k+2\right)  t$ vertices, with $\delta\left(  G\right)  =\left(  k+1\right)
n/\left(  3k+2\right)  ,$ which is not homomorphic to $A_{k}.$ Note also that
all graphs satisfying the premises of Theorem\ A are $3$-chromatic. Addressing
this last issue, Jin \cite{Jin95}, and Chen, Jin and Koh \cite{CJK97} gave a
finer characterization of $K_{3}$-free graphs with\emph{ }$\delta
>n/3.$\medskip

\textbf{Theorem B }\emph{Let }$G$\emph{ be a }$K_{3}$\emph{-}$f\emph{ree}$
\emph{graph of order }$n,$ \emph{with }$\delta\left(  G\right)  >n/3.$\emph{
If }$\chi\left(  G\right)  \geq4,$\emph{ then }$M_{3}\subset G.$ \emph{If
}$\chi\left(  G\right)  =3$\emph{ and }%
\[
\delta\left(  G\right)  >\frac{k+1}{3k+2}n,
\]
\emph{then }$G$\emph{ is homomorphic to }$A_{k}.$\medskip

Later Brandt and Pisanski \cite{BrPi98} found an infinite family of
$4$-chromatic, $K_{3}$-free graphs with $\delta\left(  G\right)  >n/3$. All
these interesting results shed some light on the structure of dense $K_{3}%
$-free graphs, but could not answer the original question of Erd\H{o}s and
Simonovits. The answer was given by Brandt and Thomass\'{e} \cite{BrTo10} in
the following ultimate result, culminating the series \cite{Bra02},
\cite{Jin95} and \cite{Tho02}.\medskip

\textbf{Theorem C }\emph{Let }$G$\emph{ be a }$K_{3}$\emph{-free graph of
order }$n.$\emph{ If }$\delta\left(  G\right)  >n/3,$\emph{ then }$\chi\left(
G\right)  \leq4.$\medskip

This theorem essentially concludes the study of $\psi\left(  n,2,h\right)  .$
Although there are still unsettled questions about $4$-chromatic $K_{3}$-free
graphs $G$ of order $n$ with $\delta\left(  G\right)  >n/3,$ the broad picture
is already fixed. The goal of this paper is to conclude likewise the study of
$\psi\left(  n,r,h\right)  $ for $r\geq3$.

To this end, in Section \ref{extHaj}, we extend the example of Hajnal and show
that%
\begin{equation}
\psi\left(  n,r,h\right)  =\left(  1-2/\left(  2r-1\right)  \right)
n-o\left(  n\right)  .\label{lob}%
\end{equation}
for every $h>r+2.$ This is to say, for every $\varepsilon$ there exists a
$K_{r+1}$-free graph of order $n$ with
\[
\delta\left(  G\right)  >\left(  1-\frac{2}{2r-1}-\varepsilon\right)  n
\]
and arbitrary large chromatic number, provided $n$ is sufficiently large. We
believe that for $r\geq3$ this extension is not widely known, although its
main idea is the same as for $r=2$.

Thus, from this point on, we are concerned mainly with the question: \emph{how
large can be }$\chi\left(  G\right)  $\emph{ of a }$K_{r+1}$\emph{-free graph
}$G$\emph{ of order }$n$\emph{ }$\emph{with}$ $\delta\left(  G\right)
>\left(  1-2/\left(  2r-1\right)  \right)  n.$ To give the reader an immediate
clue we first state an extension of Theorem C.

\begin{theorem}
\label{extBT} Let $r\geq2$ and $G$ be a $K_{r+1}$-free graph of order $n.$ If
\[
\delta\left(  G\right)  >\left(  1-\frac{2}{2r-1}\right)  n,
\]
then $\chi\left(  G\right)  \leq r+2.$
\end{theorem}

This theorem leaves only two cases of $\chi\left(  G\right)  $ to investigate,
viz., $\chi\left(  G\right)  =r+1$ and $\chi\left(  G\right)  =r+2$. As one
can expect, when $\delta\left(  G\right)  $ is sufficiently large, we have
$\chi\left(  G\right)  =r+1.$ The precise statement extends Theorem A as follows.

\begin{theorem}
\label{extJin}Let $r\geq2,$ $1\leq k\leq9,$ and let $G$ be a $K_{r+1}$-free
graph of order $n.$ If
\[
\delta\left(  G\right)  >\left(  1-\frac{2k-1}{\left(  2k-1\right)
r-k+1}\right)  n
\]
then $G$ is homomorphic to $A_{k}+K_{r-2}.$
\end{theorem}

As a corollary, under the premises of Theorem \ref{extJin}, we find that
$\chi\left(  G\right)  \leq r+1.$ Also Theorem \ref{extJin} is best possible
in the following sense: for every $k$ and $n$ there exists a $\left(
r+1\right)  $-chromatic $K_{r+1}$-free $G$ of order $n$ with%
\[
\delta\left(  G\right)  \geq\left(  1-\frac{2k-1}{\left(  2k-1\right)
r-k+1}\right)  n-1
\]
that is not homomorphic to $A_{k}+K_{r-2}.$ This example is given in Section
\ref{extAnd}.

We also generalize the example of H\"{a}ggkvist, constructing for every $n$ an
$\left(  r+2\right)  $-chromatic, $K_{r+1}$-free graph $G$ with
\[
\delta\left(  G\right)  \geq\left(  1-\frac{19}{19r-9}\right)  n-1,
\]
which shows that the conclusion of Theorem \ref{extJin} does not necessarily
hold for $k\geq10.$ This example is given in Section \ref{extHag}.

To give some further structural information, we extend Theorem C as follows.

\begin{theorem}
\label{extCJK}Let $r\geq2$ and $G$ be a $K_{r+1}$-free graph of order $n$
with
\[
\delta\left(  G\right)  >\left(  1-\frac{2}{2r-1}\right)  n.
\]
If $\chi\left(  G\right)  \geq r+2,$ then $M_{3}+K_{r-2}\subset G.$ If
$\chi\left(  G\right)  \leq r+1$ and%
\[
\delta\left(  G\right)  >\left(  1-\frac{2k-1}{\left(  2k-1\right)
r-k+1}\right)  n
\]
then $G$ is homomorphic to $A_{k}+K_{r-2}.$
\end{theorem}

This result is best possible in view of the examples described prior to
Theorem \ref{extCJK}.

For $r\geq3$ we obtain the following summary for $\psi\left(  n,r,h\right)
:$
\begin{align*}
\left(  1-3/\left(  3r-1\right)  \right)  n-1  &  \leq\psi\left(
n,r,r+1\right)  \leq\left(  1-3/\left(  3r-1\right)  \right)  n\\
\left(  1-19/\left(  19r-9\right)  \right)  n-1  &  \leq\psi\left(
n,r,r+2\right)  \leq\left(  1-19/\left(  19r-9\right)  \right)  n\\
\psi\left(  n,r,h\right)   &  =\left(  1-2/\left(  2r-1\right)  \right)
n-o\left(  n\right)  \text{ for all }h>r+2.
\end{align*}

\subsubsection*{About the proof method}

We deduce the proofs of Theorems \ref{extBT}, \ref{extJin} and \ref{extCJK} by
induction on $r$ from Theorems C, A and B respectively. Although this method
seems simple and natural, to our best knowledge none of Theorems \ref{extBT},
\ref{extJin} and \ref{extCJK} has been mentioned in the existing literature.
The same is true for the extensions of the examples of Hajnal, H\"{a}ggkvist
and Andr\'{a}sfai.

The induction step, carried uniformly in all the three proofs, is based on the
crucial Lemma \ref{reduL}. This lemma can be applied immediately to extend
other results about triangle-free graphs, but we leave these extensions to the
interested reader.

\section{Proofs}

For a graph $G$ and a vertex $u\in V\left(  G\right)  $ we write
$\Gamma\left(  u\right)  $ for the set of neighbors of $u,$ and $d_{G}\left(
u\right)  $ for $\left\vert \Gamma\left(  u\right)  \right\vert .$ If $U$ is a
subgraph of $G,$ we set- $\Gamma\left(  U\right)  =\cap_{x\in U}\Gamma\left(
x\right)  $ and $d\left(  U\right)  =\left\vert \Gamma\left(  U\right)
\right\vert ;$ note that this is not the usual definition.

Our main proof device is the following lemma.

\begin{lemma}
\label{reduL}Let $r\geq3$ and $G$ be a maximal $K_{r+1}$-free graph of order
$n.$ If
\[
\delta\left(  G\right)  >\left(  1-\frac{2}{2r-1}\right)  n,
\]
then $G$ has a vertex $u$ such that $e\left(  G_{u}^{\prime}\right)  =0$.
\end{lemma}

\textbf{Proof }For short, set $\delta=\delta\left(  G\right)  $ and
$V=V\left(  G\right)  $. Write $\delta^{\prime}\left(  H\right)  $ for the
minimum nonzero vertex degree in a graph $H,$ and set $\delta^{\prime}%
=\min\left\{  \delta^{\prime}\left(  G_{u}\right)  :u\in V\right\}  .$

We start with two facts which will be used several times throughout the proof.
First, for every set $S$ of $s\leq r-1$ vertices, we have%
\begin{align*}
d\left(  S\right)   &  \geq\sum_{v\in S}d\left(  v\right)  -\left(
s-1\right)  n\geq s\delta-\left(  s-1\right)  n\geq\left(  r-1\right)
\delta-\left(  r-2\right)  n\\
&  \geq\left(  r-1\right)  \left(  1-2/\left(  2r-1\right)  \right)  n-\left(
r-2\right)  n=n/\left(  2r-1\right)  >0;
\end{align*}
thus, every clique is contained in an $r$-clique.

Second, since $G$ is a maximal $K_{r+1}$-free graph, there is an $\left(
r-1\right)  $-clique in the common neighborhood $\Gamma\left(  uw\right)  $ of
every two nonadjacent vertices $u$ and $w$.

Now, for a contradiction, assume the conclusion of the lemma false: let
$E\left(  G_{u}^{\prime}\right)  \neq\varnothing$ for every $u\in V\left(
G\right)  $. For convenience the first part of the proof is split into
separate claims.

\begin{claim}
\label{cl1}For every edge $uv\in E\left(  G\right)  $ we have $d\left(
uv\right)  \leq\left(  r-2\right)  \left(  n-\delta\right)  -\delta^{\prime}.$
\end{claim}

\begin{proof}
Let $uv\in E\left(  G\right)  ,$ select an $\left(  r-1\right)  $-clique $R$
containing $uv,$ and let $w\in\Gamma\left(  R\right)  .$ Since $G$ is
$K_{r+1}$-free, $\Gamma\left(  R\right)  $ is an independent set, and so
$\Gamma\left(  R\right)  \subset G_{w}^{\prime}.$ It is easy to see that
$G_{w}^{\prime}$ contains at least $\delta^{\prime}$ vertices which do not
belong to $\Gamma\left(  R\right)  .$ Indeed, by assumption $G_{w}^{\prime}$
contains edges. If some of this edges contains a vertex $t\in\Gamma\left(
R\right)  ,$ we have $\Gamma_{G_{w}^{\prime}}\left(  t\right)  \cap
\Gamma\left(  R\right)  =\varnothing,$ and the assertion follows since
$d_{G_{w}^{\prime}}\left(  t\right)  \geq\delta^{\prime}.$ If no edge of
$G_{w}^{\prime}$ is incident to $\Gamma\left(  R\right)  ,$ there are at least
$\delta^{\prime}+1$ vertices of $G_{w}^{\prime}$ which do not belong to
$\Gamma\left(  R\right)  .$ Hence,%
\[
d\left(  R\right)  \leq n-d\left(  w\right)  -\delta^{\prime}\leq
n-\delta-\delta^{\prime}.
\]
Letting $S=R-u-v,$ we have
\[
d\left(  S\right)  \geq\left(  r-3\right)  \delta-\left(  r-4\right)  n,
\]
and so%
\[
n-\delta-\delta^{\prime}\geq d\left(  R\right)  \geq d\left(  uv\right)
+d\left(  S\right)  -n\geq d\left(  uv\right)  -\left(  r-3\right)  \left(
n-\delta\right)  ,
\]
as claimed.
\end{proof}

\begin{claim}
\label{cl2}For every $u\in V$ and $w\in V(G_{u}^{\prime})$ we have
\[
d\left(  uw\right)  \geq d\left(  u\right)  +d\left(  w\right)  -\left(
r-1\right)  \left(  n-\delta\right)  .
\]

\end{claim}

\begin{proof}
Since $G$ is a maximal $K_{r+1}$-free graph, we can select an $\left(
r-1\right)  $-clique $R\subset\Gamma\left(  uw\right)  $. Since $G$ is
$K_{r+1}$-free we have $\Gamma\left(  u\right)  \cap$ $\Gamma\left(  R\right)
=\varnothing$ and $\Gamma\left(  w\right)  \cap$ $\Gamma\left(  R\right)
=\varnothing$, implying that $\Gamma\left(  u\right)  \cup\Gamma\left(
w\right)  \subset V\backslash\Gamma\left(  R\right)  $, and so,
\begin{align*}
d\left(  uw\right)   &  =\left\vert \Gamma\left(  u\right)  \cap\Gamma\left(
w\right)  \right\vert =\left\vert \Gamma\left(  u\right)  \right\vert
+\left\vert \Gamma\left(  w\right)  \right\vert -\left\vert \Gamma\left(
u\right)  \cup\Gamma\left(  w\right)  \right\vert \\
&  \geq d\left(  u\right)  +d\left(  w\right)  -\left\vert V\backslash
\Gamma\left(  R\right)  \right\vert =d\left(  u\right)  +d\left(  w\right)
-n+d\left(  R\right) \\
&  \geq d\left(  u\right)  +d\left(  w\right)  -\left(  r-1\right)  \left(
n-\delta\right)  ,
\end{align*}
completing the proof.
\end{proof}

\begin{claim}
\label{cl3}For every vertex $u\in V$ the graph $G_{u}^{\prime}$ is triangle-free.
\end{claim}

\begin{proof}
Assume the opposite: let $u\in V$ and $T$ be a triangle in $G_{u}^{\prime}$
with vertices $v_{1},v_{2},v_{3}$. Using Claim \ref{cl2}, we have
\begin{align*}
d\left(  T\right)   &  =\left\vert \Gamma\left(  v_{1}\right)  \cap
\Gamma\left(  v_{2}\right)  \cap\Gamma\left(  v_{3}\right)  \right\vert
\geq\left\vert \Gamma\left(  v_{1}\right)  \cap\Gamma\left(  v_{2}\right)
\cap\Gamma\left(  v_{3}\right)  \cap\Gamma\left(  u\right)  \right\vert \\
&  \geq d\left(  v_{1}u\right)  +d\left(  v_{2}u\right)  +d\left(
v_{3}u\right)  -2d\left(  u\right) \\
&  \geq3\left(  d\left(  u\right)  -\left(  r-1\right)  \left(  n-\delta
\right)  \right)  +d\left(  v_{1}\right)  +d\left(  v_{2}\right)  +d\left(
v_{3}\right)  -2d\left(  u\right) \\
&  \geq d\left(  u\right)  +3r\delta-3\left(  r-1\right)  n\\
&  \geq\left(  3r+1\right)  \delta-3\left(  r-1\right)  n.
\end{align*}

Select an $r$-clique $R$ containing $T,$ and let $S=R-v_{1}-v_{2}-v_{3}.$ We
have
\[
d\left(  S\right)  \geq\left(  r-3\right)  \delta-\left(  r-4\right)  n,
\]
and hence,%
\begin{align*}
d\left(  R\right)   &  \geq d\left(  T\right)  +d\left(  S\right)  -n>\left(
3r+1\right)  \delta-3\left(  r-1\right)  n+\left(  r-3\right)  \delta-\left(
r-3\right)  n\\
&  \geq\left(  4r-2\right)  \delta-\left(  4r-6\right)  n>\left(  4r-2\right)
\left(  \frac{2r-3}{2r-1}\right)  n-\left(  4r-6\right)  n=0.
\end{align*}
Therefore $d\left(  R\right)  >0$, and so, $K_{r+1}\subset G,$ a contradiction
completing the proof of the claim.
\end{proof}

Let $u\in V$ and $vw\in E\left(  G_{u}^{\prime}\right)  $ be such that
$\delta^{\prime}=d_{G_{u}^{\prime}}\left(  v\right)  .$ Since $G_{u}^{\prime}$
is triangle-free, the set $\left(  \Gamma\left(  v\right)  \cap\Gamma\left(
w\right)  \right)  \backslash\Gamma\left(  u\right)  $ is empty, and so
\begin{align*}
d_{G_{u}^{\prime}}\left(  v\right)  +d_{G_{u}^{\prime}}\left(  w\right)   &
=\left\vert \Gamma\left(  v\right)  \backslash\Gamma\left(  u\right)
\right\vert +\left\vert \Gamma\left(  w\right)  \backslash\Gamma\left(
u\right)  \right\vert =\left\vert \left(  \Gamma\left(  u\right)  \cup
\Gamma\left(  w\right)  \right)  \backslash\Gamma\left(  u\right)  \right\vert
\\
&  \geq\left\vert \Gamma\left(  u\right)  \cup\Gamma\left(  w\right)
\right\vert -d\left(  u\right)  =d\left(  v\right)  +d\left(  w\right)
-\left\vert \Gamma\left(  v\right)  \cap\Gamma\left(  w\right)  \right\vert
-d\left(  u\right) \\
&  \geq2\delta-\left\vert \Gamma\left(  v\right)  \cap\Gamma\left(  w\right)
\right\vert -d\left(  u\right)  .
\end{align*}
Now,estimating $\left\vert \Gamma\left(  v\right)  \cap\Gamma\left(  w\right)
\right\vert $ by Claim \ref{cl1}, we find that
\[
\delta^{\prime}+d_{G_{u}^{\prime}}\left(  w\right)  \geq2\delta-\left(
r-2\right)  \left(  n-\delta\right)  +\delta^{\prime}-d\left(  u\right)  ,
\]
and so,
\[
d_{G_{u}^{\prime}}\left(  w\right)  \geq r\delta-\left(  r-2\right)
n-d\left(  u\right)  .
\]
On the other hand, using Claim \ref{cl2} to estimate $d\left(  uw\right)  $,
we see that%
\begin{align*}
d_{G_{u}^{\prime}}\left(  w\right)   &  =\left\vert \Gamma\left(  w\right)
\backslash\Gamma\left(  u\right)  \right\vert =d\left(  w\right)  -d\left(
uw\right) \\
&  \leq d\left(  w\right)  -d\left(  u\right)  -d\left(  w\right)  +\left(
r-1\right)  \left(  n-\delta\right) \\
&  =-d\left(  u\right)  +\left(  r-1\right)  \left(  n-\delta\right)  .
\end{align*}
Therefore,
\[
-d\left(  u\right)  +\left(  r-1\right)  \left(  n-\delta\right)  \geq
d_{G_{u}^{\prime}}\left(  w\right)  \geq r\delta-\left(  r-2\right)
n-d\left(  u\right)  ,
\]
and so
\[
\left(  2r-3\right)  n\geq\left(  2r-1\right)  \delta,
\]
a contradiction, completing the proof of the lemma.$\hfill\square$\bigskip

\begin{proof}
[\textbf{Proof of Theorem \ref{extBT}}]We shall show that $G=H+G_{0}$, where
$G_{0}$ is an $\left(  r-2\right)  $-partite graph, and $H$ is a $K_{3}$-free
graph with $\delta\left(  H\right)  >\left\vert H\right\vert /3.$ We shall
prove this assertion by induction on $r$. For $r=2$ there is nothing to prove,
so assume that the assertion holds for $r^{\prime}<r.$ Add some edges to make
$G$ maximal $K_{r+1}$-free; $\delta\left(  G\right)  $ can only increase, and
$G$ remains $K_{r+1}$-free. Lemma \ref{reduL} implies that there is a vertex
$u\in V\left(  G\right)  $ such that $G_{u}$ is empty. This means that $G$ is
homomorphic to $G_{u}+N,$ where $N$ is the graph induced by the neighbors of
$u,$ which is obviously $K_{r}$-free. We also have%
\begin{align*}
\delta\left(  N\right)   &  \geq\delta+d\left(  u\right)  -n>d\left(
u\right)  -\frac{2}{2r-1}n>d\left(  u\right)  -\frac{2}{2r-3}\delta\\
&  \geq\left(  1-\frac{2}{2r-3}\right)  \left\vert N\right\vert .
\end{align*}
By the induction hypothesis, $N$ is a join of an $\left(  r-3\right)
$-partite graph $N_{0}$, and a $K_{3}$-free graph $H$ with $\delta\left(
H\right)  >\left\vert H\right\vert /3.$ Thus $G=H+\left(  N_{0}+G_{u}\right)
,$ completing the induction step and the proof of the assertion. Since
$\chi\left(  H\right)  \leq4$, it follows that $\chi\left(  G\right)  \leq
\chi\left(  H\right)  +r-2\leq r+2,$ completing the proof.
\end{proof}

\begin{proof}
[\textbf{Proof of Theorem \ref{extJin}}]In this case we shall show that
$G=H+G_{0}$, where $G_{0}$ is an $\left(  r-2\right)  $-partite graph, and $H$
is a $K_{3}$-free graph with%
\[
\delta\left(  H\right)  >\frac{k+1}{2k+1}\left\vert H\right\vert .
\]
This assertion follows as in the proof of Theorem \ref{extBT}. The only
difference is given in the following calculation%
\begin{align*}
\delta\left(  N\right)   &  \geq\delta+d\left(  u\right)  -n>d\left(
u\right)  -\frac{2k-1}{\left(  2k-1\right)  r-k+1}n\\
&  \geq d\left(  u\right)  -\frac{1}{\left(  2k-1\right)  \left(  r-1\right)
-k+1}\delta\\
&  \geq\left(  1-\frac{1}{\left(  2k-1\right)  \left(  r-1\right)
-k+1}\right)  \left\vert N\right\vert .
\end{align*}
According to Theorem A, $H$ is homomorphic to $A_{k},$ and so $G$ is
homomorphic to $A_{k}+K_{r-2}$, completing the proof.
\end{proof}

The proof of Theorem \ref{extCJK} is the same as of Theorem \ref{extJin}, so
we shall omit it.

\section{Extension of some basic examples}

Below we construct three types of graphs by the same simple method: we take
the join of a known $K_{3}$-free graph and the $\left(  r-2\right)  $-partite
Tur\'{a}n graph. Choosing appropriately the order of the two graphs, the
resulting graph can be made almost regular.

\subsection{\label{extHaj}Extending the example of Hajnal}

In this section we shall construct, for every $r\geq2$, $h>1,$ $\varepsilon>0$
and $n$ sufficiently large, a $K_{r+1}$-free graph $G$ of order $n$ with
\begin{equation}
\delta\left(  G\right)  >\left(  1-\frac{2}{2r-1}-\varepsilon\right)  n
\label{minbo}%
\end{equation}
and $\chi\left(  G\right)  >h.$

We start by an updated version of the example of Hajnal, reported in
\cite{ErSi73}: a $K_{3}$-free graph $G$ of order $n$ with arbitrary large
chromatic number and $\delta\left(  G\right)  >n/3-o\left(  n\right)  $.

Let $K_{2m+h}^{\left(  m\right)  }$ be a Kneser graph: its vertices are the
sets $S\subset\left[  2m+h\right]  $ of size $\left\vert S\right\vert =m;$ two
vertices $S_{1}$ and $S_{2}$ are joined if $S_{1}\cap S_{2}=\varnothing.$
Clearly if $m>h$, the graph $K_{2m+h}^{\left(  m\right)  }$ is $K_{3}$-free.
Kneser \cite{Kne58} conjectured and Lov\'{a}sz \cite{Lov78} proved that
$\chi\left(  K_{2m+h}^{\left(  m\right)  }\right)  =h+2.$

\begin{example}
\label{exEHS} Let $K$ be a copy of $K_{2m+h}^{\left(  m\right)  }$ and let
$S_{1},...,S_{t}$ be its vertices, where $t=\binom{2m+h}{m}.$ Let
$n\geq3m+h+\binom{2m+h}{m},$ and set%
\[
n_{1}=n-\binom{2m+h}{m}\text{ and }k=\left\lfloor \frac{n_{1}}{3m+h}%
\right\rfloor .
\]
Add additional $n_{1}$ vertices to $K$ in the following way: add a set $A$ of
$\left(  2m+h\right)  k$ vertices, indexed for convenience as $v_{ij},$
$i\in\left[  2m+h\right]  ,$ $j\in\left[  k\right]  ,$ and add a set $B$ of
additional $n_{1}-\left(  2m+h\right)  k$ vertices. Now join every vertex of
$A$ to every vertex of $B,$ and join every vertex $v_{ij}\in A$ to every
vertex $S_{l}$ such that $i\in S_{l}.$ Write $H\left(  n,m,h\right)  $ for the
resulting graph.
\end{example}

We immediately see that $v\left(  H\left(  n,m,h\right)  \right)  =n$ and that%
\[
\chi\left(  H\left(  n,m,h\right)  \right)  \geq\chi\left(  K_{2m+h}^{\left(
m\right)  }\right)  =h+2.
\]

Let us check that $H\left(  n,m,h\right)  $ is $K_{3}$-free. Since no vertex
in $B$ is connected to a vertex in $K,$ and $A$ and $B$ are independent, after
a brief inspection, we see that a triangle in $H\left(  n,m,h\right)  $ must
have an edge $S_{i}S_{j}$ in $K$ and a vertex $v_{pq}\in A;$ thus $p\in S_{i}$
and $p\in S_{j},$ and so $S_{i}\cap S_{j}\neq\varnothing,$ contrary to the
assumption that $S_{i}S_{j}$ is an edge in $K.$ Hence, $H\left(  n,m,h\right)
$ is $K_{3}$-free.

To estimate $\delta\left(  H\left(  n,m,h\right)  \right)  $ observe that
every set $S_{i}\in V\left(  K\right)  $ is joined to $mk$ vertices of $A;$
every vertex from $B$ is joined to $\left(  2m+h\right)  k$ vertices of $A$
and every vertex of $A$ is joined to $n_{1}-\left(  2m+h\right)  k\geq mk$
vertices of $B.$ Therefore, selecting $m$ sufficiently large with respect to
$h$, we see that%
\[
\delta\left(  H\left(  n,m,h\right)  \right)  \geq mk=m\left\lfloor \left(
n-\binom{2m+h}{m}\right)  /\left(  3m+h\right)  \right\rfloor =n/3+o\left(
n\right)  .
\]
Therefore $H\left(  n,m,h\right)  $ has the required properties.

For $r\geq3$ we construct our graph $G$ as a join of a properly selected graph
$H\left(  n^{\prime},m^{\prime},h^{\prime}\right)  $ and an $\left(
r-2\right)  $-partite Tur\'{a}n graph.

\begin{example}
\label{exEHS2} Let $h>r;$ select $n_{0}$ and $m$ such that, for $n_{1}\geq
n_{0},$ we have
\[
\delta\left(  H\left(  n_{1},m,h-r\right)  \right)  >n_{1}\left(  \frac{1}%
{3}-\frac{\varepsilon}{3}\right)  .
\]
Assume that%
\[
n>\frac{2r-1}{3}n_{0};
\]
set
\begin{align*}
G_{1} &  =T_{r-2}\left(  \left\lfloor \frac{2r-4}{2r-1}n\right\rfloor \right)
\\
G_{2} &  =H\left(  n-\left\lfloor \frac{2r-4}{2r-1}n\right\rfloor
,m,h-r\right)  ,
\end{align*}
and let $G=G_{1}+G_{2}.$
\end{example}

Let us show that $G$ satisfies the requirements. Since $G_{1}$ is $K_{r-1}%
$-free and $G_{2}$ is triangle-free, we see that $G$ is $K_{r+1}$-free. Also,
we have%
\[
\chi\left(  G\right)  \geq\chi\left(  G_{2}\right)  +r-2\geq h.
\]
For every $v\in G_{1},$
\begin{align*}
d\left(  v\right)   &  =d_{G_{1}}\left(  v\right)  +\left\vert G_{2}%
\right\vert \geq\left\lfloor \frac{r-3}{r-2}\left\lfloor \frac{2r-4}%
{2r-1}n\right\rfloor \right\rfloor +n-\left\lfloor \frac{2r-4}{2r-1}%
n\right\rfloor \\
&  \geq n-\frac{2}{2r-1}n-1.
\end{align*}
On the other hand, for every $v\in G_{2}$,%
\begin{align*}
d\left(  v\right)   &  =d_{G_{2}}\left(  v\right)  +\left\vert G_{1}%
\right\vert \geq\left\lfloor \frac{2r-4}{2r-1}n\right\rfloor +\left(
n-\left\lfloor \frac{2r-4}{2r-1}n\right\rfloor \right)  \left(  \frac{1}%
{3}-\varepsilon\right)  \\
&  \geq\frac{2r-4}{2r-1}n-1+\frac{3r}{2r-1}n\left(  \frac{1}{3}-\frac
{\varepsilon}{3}\right)  =\left(  1-\frac{2}{2r-1}-\frac{r\varepsilon}%
{2r-1}-\frac{1}{n}\right)  n\\
&  >\left(  1-\frac{2}{2r-1}-\varepsilon\right)  n.
\end{align*}
Hence, (\ref{minbo}) also holds, and thus, $G$ has the required properties.

From our construction and Theorem \ref{extBT} it follows that for all $h>r+2,$%
\[
\psi\left(  n,r,h\right)  =\left(  1-2/\left(  2r-1\right)  \right)
n-o\left(  n\right)  .
\]

\subsection{\label{extHag}Extending the example of H\"{a}ggkvist}

As mentioned in the introduction, H\"{a}ggkvist\cite{Hag82} constructed for
every $k\geq1,$ a $4$-chromatic, $10k$-regular graph of order $29k.$ For
completeness we describe this example.

\begin{example}
\label{exHag}Partition $V\left(  G\right)  =\left[  n\right]  $ into $11$
sets
\[
\left[  n\right]  =A_{1}\cup...\cup A_{5}\cup B_{1}\cup...\cup B_{5}\cup C
\]
such that $\left\vert A_{1}\right\vert =\cdots=\left\vert A_{5}\right\vert
=3k,$ $\left\vert B_{1}\right\vert =\cdots=\left\vert B_{5}\right\vert =2k,$
$\left\vert C\right\vert =4k;$ join $u\in A_{i}$ to $v\in A_{j}$ if $i-j=\pm1$
$\operatorname{mod}$ $5$; join $u\in A_{i}$ to $v\in B_{j}$ if $i-j=\pm1$
$\operatorname{mod}$ $5;$ join all vertices of $C$ to all vertices of
$\cup_{i=1}^{5}B_{i}.$

Write $H\left(  k\right)  $ for the resulting graph.
\end{example}

Observe that $H\left(  k\right)  $ contains the Mycielski graph $M_{3},$ which
is $K_{3}$-free and $4$-chromatic$.$ In fact, $H\left(  k\right)  $ is
homomorphic to $M_{3};$ hence, it is $K_{3}$-free and $4$-chromatic itself. It
is obvious that $\delta\left(  H\left(  k\right)  \right)  =10k.$

Now we shall construct for every $r\geq3$ and every $n>19r-9$ a $K_{r+1}%
$-free, $\left(  r+2\right)  $-chromatic graph of order$\ n$ with%
\begin{equation}
\delta>\left(  1-\frac{19}{19r-9}\right)  n-1. \label{minbo2}%
\end{equation}

\begin{example}
Assume that $n>19r-9;$ set%
\begin{align*}
G_{1}  &  =T_{r-2}\left(  n-29\left\lfloor \frac{n}{19r-9}\right\rfloor
\right) \\
G_{2}  &  =A_{k}\left(  \left\lfloor \frac{n}{19r-9}\right\rfloor \right)  ,
\end{align*}
and let $G=G_{1}+G_{2}.$
\end{example}

We shall show that $G$ satisfies the requirements. Since $G_{1}$ is $K_{r-1}%
$-free and $G_{2}$ is $K_{3}$-free, we see that $G$ is $K_{r+1}$-free. Also,
we have%
\[
\chi\left(  G\right)  =\chi\left(  G_{2}\right)  +r-2=r+2.
\]
For every $v\in G_{2},$%
\begin{align*}
d\left(  v\right)   &  =d_{G_{2}}\left(  v\right)  +\left\vert G_{1}%
\right\vert \geq n-29\left\lfloor \frac{n}{19r-9}\right\rfloor +10\left\lfloor
\frac{n}{19r-9}\right\rfloor \\
&  \geq\left(  1-\frac{19}{19r-9}\right)  n.
\end{align*}
On the other hand, for every $v\in G_{1}$ we have
\begin{align*}
d\left(  v\right)   &  =d_{G_{1}}\left(  v\right)  +\left\vert G_{2}%
\right\vert \geq\delta\left(  T_{r-2}\left(  n-29\left\lfloor \frac{n}%
{19r-9}\right\rfloor \right)  \right)  +29\left\lfloor \frac{n}{19r-9}%
\right\rfloor \\
&  =\left\lfloor \frac{r-3}{r-2}\left(  n-29\left\lfloor \frac{n}%
{19r-9}\right\rfloor \right)  \right\rfloor +29\left\lfloor \frac{n}%
{19r-9}\right\rfloor \\
&  =\left\lfloor \frac{r-3}{r-2}n+\frac{29}{r-2}\left\lfloor \frac{n}%
{19r-9}\right\rfloor \right\rfloor .
\end{align*}
Suppose that $n=\left(  19r-9\right)  k+s,$ where $k\geq0$ and $0\leq
s\leq19r-10$ are integers. Then
\begin{align*}
\left\lfloor \frac{r-3}{r-2}n+\frac{29}{r-2}\left\lfloor \frac{n}%
{19r-9}\right\rfloor \right\rfloor  &  =n+\left\lfloor -\frac{\left(
19-9\right)  k+s}{r-2}+\frac{29k}{r-2}\right\rfloor \\
&  =n-19k+\left\lfloor -\frac{s}{r-2}\right\rfloor >n-19k-\frac{s}{r-2}-1\\
&  =\left(  1-\frac{19}{19r-9}\right)  n-1.
\end{align*}

Hence, (\ref{minbo2}) also holds, and $G$ has the required properties. Note
that if $19r-9$ divides $n,$ then
\[
\delta\left(  G\right)  =\left(  1-\frac{19}{19r-9}\right)  n.
\]

\subsection{\label{extAnd}Extending the Andr\'{a}sfai graphs}

Let $A_{k}$ be the $k$'th Andr\'{a}sfai graph, which is a $k$-regular graph of
order $3k-1.$ Write $A_{k}\left(  t\right)  ,$ for the blow-up of $A_{k}$ by
factor $t,$ i.e., $A_{k}\left(  t\right)  $ is obtained by replacing each
vertex $u\in V\left(  A_{k}\right)  $ with a set $V_{u}$ of size $t$ and each
edge $uv\in E\left(  H\right)  $ with a complete bipartite graph with vertex
classes $V_{u}$ and $V_{v}.$\ Note that $A_{k}\left(  t\right)  $ is $K_{3}%
$-free, $3$-chromatic $kt$-regular graph of order $\left(  3k-1\right)  t$.

We shall construct for every $r\geq3$ and every $n>\left(  2k-1\right)  r-k+1$
a $K_{r+1}$-free, $\left(  r+1\right)  $-chromatic graph of order$\ n$ with%
\begin{equation}
\delta>\left(  1-\frac{2k-1}{\left(  2k-1\right)  r-k+1}\right)
n-1.\label{minbo3}%
\end{equation}

\begin{example}
Assume that $n>\left(  2k-1\right)  r-k+1;$ set%
\begin{align*}
G_{1} &  =T_{r-2}\left(  n-\left(  3k-1\right)  \left\lfloor \frac{n}{\left(
2k-1\right)  r-k+1}\right\rfloor \right)  \\
G_{2} &  =A_{k}\left(  \left\lfloor \frac{n}{\left(  2k-1\right)
r-k+1}\right\rfloor \right)  ,
\end{align*}
and let $G=G_{1}+G_{2}.$
\end{example}

We shall show that $G$ satisfies the requirements. Since $G_{1}$ is $K_{r-1}%
$-free and $G_{2}$ is $K_{3}$-free, we see that $G$ is $K_{r+1}$-free. Also,
we have%
\[
\chi\left(  G\right)  =\chi\left(  G_{2}\right)  +r-2=r+1.
\]
For every $v\in G_{2},$%
\begin{align*}
d\left(  v\right)   &  =\left\vert G_{1}\right\vert +d_{G_{2}}\left(
v\right)  =n-\left(  3k-1\right)  \left\lfloor \frac{n}{\left(  2k-1\right)
r-k+1}\right\rfloor +k\left\lfloor \frac{n}{\left(  2k-1\right)
r-k+1}\right\rfloor \\
&  \geq\left(  1-\frac{2k-1}{\left(  2k-1\right)  r-k+1}\right)  n.
\end{align*}
On the other hand, for every $v\in G_{1}$ we have
\begin{align*}
d\left(  v\right)   &  =d_{G_{1}}\left(  v\right)  +\left\vert G_{2}%
\right\vert \\
&  \geq\delta\left(  T_{r-2}\left(  n-\left(  3k-1\right)  \left\lfloor
\frac{n}{\left(  2k-1\right)  r-k+1}\right\rfloor \right)  \right)  +\left(
3k-1\right)  \left\lfloor \frac{n}{\left(  2k-1\right)  r-k+1}\right\rfloor \\
&  =\left\lfloor \frac{r-3}{r-2}\left(  n-\left(  3k-1\right)  \left\lfloor
\frac{n}{\left(  2k-1\right)  r-k+1}\right\rfloor \right)  \right\rfloor
+\left(  3k-1\right)  \left\lfloor \frac{n}{\left(  2k-1\right)
r-k+1}\right\rfloor \\
&  \geq\left\lfloor \frac{r-3}{r-2}n+\frac{3k-1}{r-2}\left\lfloor \frac
{n}{\left(  2k-1\right)  r-k+1}\right\rfloor \right\rfloor .
\end{align*}
Suppose that $n=\left(  \left(  2k-1\right)  r-k+1\right)  t+s,$ where
$t\geq1$ and $0\leq s<\left(  2k-1\right)  r-k+1$ are integers. Then
\begin{align*}
\left\lfloor \frac{r-3}{r-2}n+\frac{3k-1}{r-2}\left\lfloor \frac{n}{\left(
2k-1\right)  r-k+1}\right\rfloor \right\rfloor  &  =n+\left\lfloor
-\frac{\left(  \left(  2k-1\right)  r-k+1\right)  t+s}{r-2}+\frac{\left(
3k-1\right)  t}{r-2}\right\rfloor \\
&  =n-\left(  2k-1\right)  t+\left\lfloor -\frac{s}{r-2}\right\rfloor \\
&  >n-\left(  2k-1\right)  t-\frac{s}{r-2}-1\\
&  =\left(  1-\frac{2k-1}{\left(  2k-1\right)  r-k+1}\right)  n-1.
\end{align*}

Hence, (\ref{minbo3}) also holds, and $G$ has the required properties. Note
that if $\left(  2k-1\right)  r-k+1$ divides $n,$ then
\[
\delta\left(  G\right)  =\left(  1-\frac{2k-1}{\left(  2k-1\right)
r-k+1}\right)  n.
\]

\textbf{Acknowledgement }This research has been supported in part by NSF Grant
\# DMS-0906634.

\end{document}